\documentclass[12pt,a4paper]{amsart}
\usepackage{amssymb,amsmath}

\headsep=1cm \numberwithin{equation}{section}
\newtheorem{theorem}{Theorem}[section]
\newtheorem{lemma}[theorem]{Lemma}
\newtheorem{proposition}[theorem]{Proposition}
\newtheorem{corollary}[theorem]{Corollary}

\theoremstyle{definition}
\newtheorem{definition}[theorem]{Definition}
\theoremstyle{remark}

\newtheorem{example}[theorem]{Example}

\newcommand{\E}{\operatorname{E}}

\newcommand{\Hom}{\operatorname{Hom}}

\newcommand{\lo}{\longrightarrow}
\newcommand{\fm}{\frak{m}}

\def\mapdown#1{\Big\downarrow\rlap
{$\vcenter{\hbox{$\scriptstyle#1$}}$}}

\begin{document}
\author[Divaani-Aazar, Esmkhani and Tousi ]{Kamran Divaani-Aazar, Mohammad Ali
Esmkhani and Massoud Tousi}
\title[Two characterizations of... ]
{Two characterizations of pure injective modules}

\address{K. Divaani-Aazar, Department of Mathematics, Az-Zahra University,
Vanak, Post Code 19834, Tehran, Iran-and-Institute for Studies in
Theoretical Physics and Mathematics, P.O. Box 19395-5746, Tehran,
Iran.} \email{kdivaani@ipm.ir}

\address{M.A. Esmkhani, Department of Mathematics, Shahid Beheshti University,
Tehran, Iran-and-Institute for Studies in Theoretical Physics and
Mathematics, P.O. Box 19395-5746, Tehran, Iran.}

\address{M. Tousi, Department of Mathematics, Shahid Beheshti University, Tehran,
Iran-and-Institute for Studies in Theoretical Physics and
Mathematics, P.O. Box 19395-5746, Tehran, Iran.}


\keywords{Pure injective modules, injective cogenerators, finitely
embedded modules, finitely presented modules.\\
This research was in part supported by a grant from IPM (No.
83130115).}

\begin{abstract}
Let $R$ be a commutative ring with identity and $D$ an $R$-module.
It is shown that if $D$ is pure injective, then $D$ is isomorphic
to a direct summand of the direct product of a family of finitely
embedded modules. As a result, it follows that if $R$ is
Noetherian, then $D$ is pure injective if and only if $D$ is
isomorphic to a direct summand of the direct product of a family
of Artinian modules. Moreover, it is proved that $D$ is pure
injective if and only if there is a family
$\{T_\lambda\}_{\lambda\in \Lambda}$ of $R$-algebras which are
finitely presented as $R$-modules, such that $D$ is isomorphic to
a direct summand of a module of the form $\Pi_{\lambda\in
\Lambda}E_\lambda$ where for each $\lambda\in \Lambda$,
$E_\lambda$ is an injective $T_\lambda$-module.
\end{abstract}

\maketitle

\section{Introduction}

Throughout this paper, let $R$ denote a commutative ring with
identity and all modules are assumed to be unitary. The notion of
injective modules has a substantial role in algebra. There are
several generalizations of this notion. One of them is the notion
of pure injective modules. An $R$-monomorphism $f:M\lo N$ is said
to be pure if for any $R$-module $L$, the map $f\otimes
id_L:M\otimes_RL\lo N\otimes_RL$ is injective. An $R$-module $D$
is said to be pure injective if for any pure homomorphism $f:M\lo
N$, the induced homomorphism $\Hom_R(N,D)\lo \Hom_R(M,D)$ is
surjective. In model theory the notion of pure injective modules
is more useful than that of injective modules. Also, there are
some excellent applications of this notion in the theory of flat
covers. Thus, this notion has attended more notice in recent
years. For a survey on pure injective modules, we refer the reader
to [{\bf 8}], [{\bf 3}] and [{\bf 9}].

Our aim in this paper is to present two characterizations of pure
injective modules. First, we show that a pure injective $R$-module
$D$ is isomorphic to a direct summand of the direct product of a
family of finitely embedded $R$-modules. Then we deduce our first
characterization that says that over a Noetherian ring $R$ an
$R$-module $D$ is pure injective if and only if $D$ is isomorphic
to a direct summand of the direct product of a family of Artinian
modules.

Let $E$ denote an injective cogenerator of the ring $R$. Let $D$
be an $R$-module. We show that $D$ is pure injective if and only
if there is a family $\{T_\lambda\}_{\lambda\in \Lambda}$ of
$R$-algebras which are finitely presented as $R$-modules, such
that $D$ is isomorphic to a direct summand of $\Pi_{\lambda\in
\Lambda}\Hom_R(T_\lambda,E)$. As an immediate consequence, it
follows that $D$ is pure injective if and only if there are an
injective $R$-module $E'$ and an $R$-module $L$, which is the
direct sum of a family of finitely presented $R$-modules such that
$D$ is isomorphic to a direct summand of $\Hom_R(L,E')$. Finally,
we will deduce our second characterization of pure injective
modules which asserts that $D$ is pure injective if and only if
there is a family $\{T_\lambda\}_{\lambda\in \Lambda}$ of
$R$-algebras which are finitely presented as $R$-modules, such
that $D$ is isomorphic to a direct summand of a module of the form
$\Pi_{\lambda\in \Lambda}E_\lambda$ where for each $\lambda\in
\Lambda$, $E_\lambda$ is an injective $T_\lambda$-module.

\section{The results}

First, we need to recall a definition and bring some lemmas.

\begin{definition} An $R$-module $M$ is called {\it cocyclic} if
$M$ is isomorphic to a submodule of the injective envelope of a
simple module (see [{\bf 5}, page 4]). We say that an $R$-module
$M$ is {\it finitely embedded} if $M$ is isomorphic to a submodule
of the injective envelope of the direct sum of finitely many
simple modules.
\end{definition}

The definition of the dual notion of ``finitely generated" are
also given separately in [{\bf 7}] and [{\bf 1}]. In the following
lemma, we show that those definitions are equivalent to the above
definition. In the sequel, for an $R$-module $M$, let $\E(M)$
denote its injective envelope.

\begin{lemma} Let $M$ be nonzero $R$-module. Then the following are
equivalent:\\
i) $M$ is finitely embedded.\\
ii) There are simple $R$-modules $S_1,S_2,\dots ,S_t$ such that
$\E(M)=\E(S_1)\oplus \dots \oplus \E(S_t)$.\\
iii) The socle of $M$ is a finitely generated and $M$ is an
essential extension of its socle.\\
iv) For any family $\{M_i\}_{i\in I}$ of submodules of $M$, the
intersection $\cap_{i\in I}M_i$ is nonzero, whenever the
intersection of any finite number of $M_i$'s is nonzero.
\end{lemma}

{\bf Proof.} $i)\Leftrightarrow ii)$ and $ii)\Leftrightarrow iii)$
follows, by [{\bf 7}, Proposition 3.20] and [{\bf 7}, Proposition
3.18] respectively. The equivalence $iii)\Leftrightarrow iv)$
follows by [{\bf 1}, Proposition 10.7]. $\Box$

\begin{lemma} Let $M$ be a nonzero $R$-module. Then
$M$ is cocyclic if and only if the intersection of all nonzero
submodules of $M$ is nonzero.
\end{lemma}

{\bf Proof.} First, assume that $M$ is cocyclic. Then, there is a
simple $R$-module $S$ such that $M$ is a submodule of $\E(S)$.
Since $\E(S)$ is an essential extension of $S$, it follows that
every nonzero submodule of $M$ contains $S$.

Conversely, assume that there is a nonzero element $c\in M$ such
that $c$ belongs to all nonzero submodules of $M$. Then, it is
routine check that $Rc$ is  simple and that $M$ is an essential
extension of $Rc$. Hence $M$ can be naturally embedded in
$\E(Rc)$. $\Box$

\begin{lemma} Let $M$ be an $R$-module. For any finitely generated
$R$-module $N$ and any $0\neq \alpha\in M\otimes_RN$, there exists
a finitely embedded quotient module $L$ of $M$ such that the image
of $\alpha$ under the natural homomorphism $M\otimes_RN\lo
L\otimes_RN$ is nonzero.
\end{lemma}

{\bf Proof.} Assume that $N$ is generated by $\{n_1,n_2,\dots
,n_s\}$. Let $F$ be a free $R$-module of rank $s$. Suppose that
$\{x_1,x_2,\dots ,x_s\}$ is a basis for $F$ and that $\rho: F\lo
N$ is the natural epimorphism defined by
$\sum_{i=1}^sr_ix_i\mapsto  \sum_{i=1}^sr_in_i$. Let $K$ denote
the kernel of $\rho$ and $\mu:K\lo F$ denote the inclusion map.
Assume $\beta=\sum_{i=1}^s(m_i\otimes x_i)\in M\otimes_RF$ is such
that $(id_M\otimes\rho)(\beta)=\alpha$.  Let $\Sigma$ denote the
set of all submodules $P$ of $M\otimes_RF$ such that
$im(id_M\otimes \mu)\subseteq P$ and $\beta\not\in P$. Since
$\beta\not\in im(id_M\otimes \mu)$, it follows that $\Sigma$ is
not empty. Clearly, $(\Sigma,\subseteq)$ satisfies the assumptions
of Zorn's Lemma, and so $\Sigma$ possesses a maximal element $Q$
'say.  Lemma 2.3  yields that $(M\otimes_RF)/Q$ is cocyclic. Fix
$1\leq i\leq s$, and let $f_i:M\lo M\otimes_RF$ be the
homomorphism defined by $m\mapsto m\otimes x_i$. Let
$f_i^*:M/f_i^{-1}(Q)\lo (M\otimes_RF)/Q$ denote the induced
monomorphism. It follows that $M/f_i^{-1}(Q)$ is cocyclic. Let
$L=M/\cap_{i=1}^sf_i^{-1}(Q)$. Since $\oplus_{i=1}^sM/f_i^{-1}(Q)$
is finitely embedded and $L$ can be naturally embedded in
$\oplus_{i=1}^sM/f_i^{-1}(Q)$, it turns out that $L$ is finitely
embedded. Let $\pi:M\lo L$ denote the natural epimorphism. The
exact sequence  $$0\lo K\overset{\mu}{\lo} F\overset{\rho}{\lo}
N\lo 0,$$ yields the following commutative diagram.

\begin{equation*}
\setcounter{MaxMatrixCols}{11}
\begin{matrix}
&M\otimes_RK &\stackrel{id_M\otimes \mu}{\lo}  &M\otimes_RF
&\stackrel{id_M\otimes \rho}{\lo} &M\otimes_RN&\lo 0
\\&\mapdown{\pi\otimes id_K} & &\mapdown{\pi\otimes id_F}
 & &\mapdown{\pi\otimes id_N} & & & &
\\  &L\otimes_RK &\stackrel{id_L\otimes \mu}{\lo}
&L\otimes_RF&\stackrel{id_L\otimes \rho}{\lo} &L\otimes_RN&\lo 0.
&
\end{matrix}
\end{equation*}

To show that $(\pi\otimes id_N)(\alpha)\neq 0$, it is enough to
see that $(\pi\otimes id_F)(\beta)\not\in im(id_L\otimes \mu)$.
Suppose the contrary is true. Then there are $y_1,y_2,\dots
,y_l\in M$ and  $z_1,z_2,\dots ,z_l\in K$ such that
$$\Sigma_{i=1}^s(\pi(m_i)\otimes x_i)=(\pi\otimes
id_F)(\beta)=(id_L\otimes \mu)(\Sigma_{j=1}^l(\pi(y_j)\otimes
z_j)).$$ For any given integer $1\leq j\leq l$, let $r_{ij}\in R,
i=1,2,\dots ,s$ be such that $z_j=\Sigma_{i=1}^sr_{ij}x_i$. Thus
$$\Sigma_{i=1}^s(\pi(m_i)\otimes
x_i)=\Sigma_{i=1}^s((\Sigma_{j=1}^l\pi(r_{ij}y_j))\otimes x_i),$$
as elements of $L\otimes_RF$. Therefore
$m_i-\Sigma_{j=1}^lr_{ij}y_j\in \bigcap_{k=1}^sf_k^{-1}(Q)$ for
all $i=1,\dots ,s$, and so
$$\Sigma_{i=1}^s(m_i\otimes
x_i)-\Sigma_{i=1}^s((\Sigma_{j=1}^lr_{ij}y_j)\otimes x_i)\in Q.$$
Because $\Sigma_{j=1}^l(y_j\otimes z_j)\in im(id_M\otimes
\mu)\subseteq Q$, it turns out that $\beta\in Q$. Therefore, we
achieved at a contradiction. $\Box$

\begin{proposition} Let $M$ be an $R$-module and let
$\{Q_i\}_{i\in I}$ denote the class of all finitely embedded
quotient modules of $M$. Then there exists a pure homomorphism
$M\lo \Pi_{i\in I}Q_i$.
\end{proposition}

{\bf Proof.}  For each $i\in I$, let $\pi_i:M\lo Q_i$ denote the
natural epimorphism. Define the homomorphism $\lambda:M\lo
\Pi_{i\in I}Q_i$, by $m\mapsto (\pi_i(m))_{i\in I}$. To show that
$\lambda$ is pure, it is enough to see that for any finitely
presented $R$-module $N$, the $R$-homomorphism
$$\lambda\otimes id_N:M\otimes_RN\lo (\Pi_{i\in I}Q_i)\otimes_RN,$$
is injective. Note that, by [{\bf 3}, Construction 18-2.6] any
$R$-module is isomorphic to the direct limit of a family of
finitely presented modules. It is a routine check to see that the
map
$$\Theta: (\Pi_{i\in I}Q_i)\otimes_RN\lo \Pi_{i\in
I}(Q_i\otimes_RN),$$ defined by $(y_i)_{i\in I}\otimes z\mapsto
(y_i\otimes z)_{i\in I}$ is an isomorphism. Now, the conclusion
follows, because by Lemma 2.4, $\Theta(\lambda\otimes id_N)$ is
injective. $\Box$

Let $D$ be a pure injective $R$-module and $f:D\lo L$ a pure
homomorphism. Then, it follows that $D$ is isomorphic to a direct
summand of $L$. This easy fact yields the following.

\begin{corollary} Let $D$ be a pure injective $R$-module. Then $D$
is isomorphic to a direct summand of the direct product of a
family of finitely embedded modules.
\end{corollary}

\begin{example} By [{\bf 5}, Theorem 6], a Pr\"{u}fer domain $R$ is locally
almost maximal if and only if every finitely embedded $R$-module
is pure injective. On other hand, by [{\bf 2}, Example 2.4] there
exists a valuation domain $R$ such that $R$ is not almost maximal.
Hence finitely embedded modules are not pure injective in general,
and so the converse of Corollary 2.6 is not true.
\end{example}

Now, we are ready to establish our first characterization of pure
injective modules.

\begin{theorem} Let $R$ be a Noetherian ring and $D$ an
$R$-module. Then $D$ is pure injective if and only if $D$ is
isomorphic to a direct summand of the direct product of a family
of Artinian modules.
\end{theorem}

{\bf Proof.} First assume that $D$ is isomorphic to a direct
summand of the direct product of a family of Artinian modules. By
[{\bf 6}, Corollary 4.2], any Artinian module is pure injective.
Thus $D$ is a pure injective module. Note that the direct product
of a family $\{D_i\}_{i\in I}$ of $R$-modules is pure injective if
and only if each $D_i$ is pure injective.

Now, assume that $D$ is pure injective. Then, it follows that $D$
is isomorphic to a direct summand of the direct product of a
family of finitely embedded $R$-modules. But, over a Noetherian
ring  finitely embedded $R$-modules are Artinian. $\Box$

\begin{example} There is a pure injective module which is not the
direct product of any family of Artinian modules. To this end, let
$(R,\fm)$ be a Noetherian complete local ring which is not
Artinian. Since $R\cong \Hom_R(\E(R/\fm),\E(R/\fm))$, it follows
by [{\bf 4}, Lemma 2.1] that $R$ is pure injective. Assume that
there is a family $\{A_i\}_{i\in I}$ of Artinian modules such that
$R=\Pi_{i\in I}A_i$. Then only finitely many of $A_i$'s can be
nonzero, because $R$ is Noetherian. Hence $R$ becomes an Artinian
$R$-module and so we achieved at a contradiction.
\end{example}

\begin{lemma} Let $\Omega$ denote the class of all $R$-algebras
$T$ which are finitely presented as $R$-modules. Let $N$ be a
submodule of an $R$-module $M$. Then the inclusion map
$i:N\hookrightarrow M$ is pure if and only if $i\otimes id_T$ is
injective for all $T\in \Omega$.
\end{lemma}

{\bf Proof.} It is enough to prove that the ``if" part. To this
end, by [{\bf 3}, Construction 18-2.6], it will suffice to show
that for any finitely presented $R$-module $L$, $i\otimes id_L$ is
injective. Let $L$ be a finitely presented $R$-module and
$T_L=R\oplus L$ denote the trivial extension of $R$ by $L$. Recall
that addition and multiplication in $T_L$ are defined respectively
by,
$$(r,x)+(r',x')=(r+r',x+x')$$
and $$(r,x)(r',x')=(rr',rx'+r'x).$$ Consider the ring homomorphism
$\psi_L:R\lo T_L$ defined by $r\mapsto (r,0)$. Clearly, the
$R$-module structure on $T_L$ that induced by $\psi_L$ coincides
with the usual $R$-module structure of $T_L$ when $T_L$ is
considered as the direct sum of two $R$-modules $R$ and $L$.

Let $\lambda:L\lo T_L$ be the natural embedding defined by
$x\mapsto (0,x)$. Now, consider the following diagram
\begin{equation*}
\setcounter{MaxMatrixCols}{11}
\begin{matrix}
&N\otimes_RL &\stackrel{i\otimes id_L}{\lo} &M\otimes_RL
\\&\mapdown{id_N\otimes \lambda}  & &\mapdown{id_M\otimes \lambda}
& & & &
\\&N\otimes_RT_L &\stackrel{i\otimes id_{T_{L}}} {\lo} &M\otimes_RT_L &
\end{matrix}
\end{equation*}
Because $\lambda$ is a pure $R$-homomorphism, we deduce that both
maps $id_N\otimes \lambda$ and $id_M\otimes \lambda$ are
injective. Thus we conclude that $i\otimes id_L$ is injective, as
required. $\Box$

\begin{proposition} Let $E$ be an injective cogenerator of $R$
and $D$ an $R$-module. Then $D$ is pure injective if and only if
there is a family $\{T_\lambda\}_{\lambda\in \Lambda}$ of
$R$-algebras which are finitely presented as $R$-modules, such
that $D$ is isomorphic to a direct summand of
$\Pi_{\lambda\in\Lambda}\Hom_R(T_\lambda,E)$.
\end{proposition}

{\bf Proof.} It is easy to see  that  for any $R$-module $M$ and
any injective $R$-module $E'$, the $R$-module $\Hom_R(M,E')$ is
pure injective (see e.g. [{\bf 4}, Lemma 2.1]). Thus the ``if"
part follows clearly, because the direct product of a family
$\{D_i\}_{i\in I}$ of $R$-modules is pure injective if and only if
$D_i$ is pure injective for all $i\in I$.

Now, we prove the converse. Let $\Omega$ be as in Lemma 2.10.
There is a subclass $\Omega^*$ of $\Omega$, which is a set , with
the property that any $T\in \Omega$ is isomorphic (as an
$R$-module) to an element of $\Omega^*$. Let $(.)^\vee$ denote the
functor $\Hom_R(.,E)$. Also, let $\Lambda$ be the set of all pairs
$(T,f)$ with $T\in \Omega^*$ and $f\in \Hom_R(D,T^\vee)$, and for
each $\lambda\in \Lambda$ denote the corresponding pair by
$(T_{\lambda},f_{\lambda})$. Set $C=\Pi_{\lambda\in
\Lambda}(T_{\lambda})^\vee$ and define the monomorphism $\psi:D\lo
C$, by $\psi(x)=(f_{\lambda}(x))_{\lambda\in \Lambda}$. By, the
paragraph preceding Corollary 2.6, the proof will be complete if
we show that $\psi$ is pure. To this end, by Lemma 2.10, it is
enough to show that for any $T\in \Omega^*$, the map $\psi\otimes
id_T$ is injective. Because, the functor $(.)^\vee$ is faithfully
exact, we equivalently prove that the natural $R$-homomorphism
$\Hom_R(C,T^\vee)\lo \Hom_R(D,T^\vee)$ is surjective.  Let $f\in
\Hom_R(D,T^\vee)$. There is $\lambda_0\in \Lambda$ with
$(T_{\lambda_0},f_{\lambda_0})=(T,f)$. Denote the projection map
$\Pi_{\lambda\in \Lambda}(T_{\lambda})^\vee\lo
(T_{\lambda_0})^\vee$, by $\rho_{\lambda_0}$. Then
$\rho_{\lambda_0}\psi=f$. $\Box$

Next, we present another characterization of pure injective
modules.

\begin{corollary} Let $D$ be an $R$-module. Then $D$ is pure
injective if and only if  $D$ is isomorphic to a direct summand of
a module of the form $\Hom_R(L,E)$ where $E$ is an injective
$R$-module and $L$ is the direct sum of a family of finitely
presented $R$-modules.
\end{corollary}

{\bf Proof.}  The ``if" part is clear, as we have seen in the
proof of Proposition 2.11.

Now, let $D$ be a pure injective $R$-module and $E$ an injective
cogenerator of $R$. Then, by Proposition 2.11, there is a family
$\{L_i\}_{i\in I}$ of finitely presented $R$-modules such that $D$
is isomorphic to a direct summand of $\Pi_{i\in I}\Hom_R(L_i,E)$.
Set $L=\oplus_{i\in I}L_i$. Then
$$\Pi_{i\in I}\Hom_R(L_i,E)\cong \Hom_R(L,E),$$and so the
conclusion follows. $\Box$

Now, we are ready to present our second characterization of pure
injective modules.

\begin{theorem}  An $R$-module $D$ is pure injective if and only
if there is a family $\{T_\lambda\}_{\lambda\in \Lambda}$ of
$R$-algebras which are finitely presented as $R$-modules, such
that $D$ is isomorphic to a direct summand of a module of the form
$\Pi_{\lambda\in \Lambda}E_\lambda$ where for each $\lambda\in
\Lambda$, $E_\lambda$ is an injective $T_\lambda$-module.
\end{theorem}

{\bf Proof.}  Again, the ``if" part is clear, because as one can
see easily, any pure injective module over an $R$-algebra $T$ is
also pure injective as an $R$-module.

Now, we prove the converse. Let $E$ be an injective $R$-module and
$T$ an $R$-algebra. Then, it is easy to see that $\Hom_R(T,E)$ is
an injective $T$-module. Hence the claim follows, by Proposition
2.11. $\Box$



\begin{thebibliography}{99}

\bibitem{} F.W. Anderson and K.R. Fuller, {\it Rings and categories of
modules}, 2nd edition, Springer-Verlag, New York, 1992.
\bibitem{} W. Brandal, {\it Almost maximal integral domains and finitely
generated modules}, Trans. Amer. Math. Soc., {\bf 183} (1973),
203-222.
\bibitem{} J. Dauns, {\it Modules and rings}, Cambridge University Press,
Cambridge, 1994.
\bibitem{} E. Enochs, {\it Flat covers and flat cotorsion modules},
Proc. Amer. Math. Soc., {\bf 92}(2) (1984), 179-184.
\bibitem{} L. Fuchs and A. Meijer, {\it Note on modules
over Pr\"{u}fer domains}, Math. Pannon, {\bf 2}(1) (1991), 3-11.
\bibitem{} L. Melkersson, {\it Cohomological properties of modules with
secondary representations}, Math. Scand., {\bf 77}(2) (1995),
197-208.
\bibitem{} D.W. Sharpe and P. V$\acute{a}$mos, {\it Injective modules},
Cambridge Tracts in Mathematics and Mathematical Physics, {\bf
62}, Cambridge University Press, London-New York, 1972.
\bibitem{} R.B. Warfield, {\it Purity and algebraic compactness for
modules}, Pacific J. Math., {\bf 28} (1969), 699-719.
\bibitem{} J. Xu, {\it Flat covers of modules}, Lecture Notes in Mathematics,
{\bf 1634}, Springer-Verlag, Berlin, 1996.

\end{thebibliography}
\end{document}